\documentclass[12pt]{article}
\usepackage{amsmath,amssymb,amsthm}
\usepackage[mathscr]{eucal}
\usepackage{rotate,graphics,epsfig}
\usepackage{color}
\usepackage{hyperref}

\newtheorem{The}{Theorem}
\newtheorem{Exa}[The]{Example}

\newtheorem{Pro}[The]{Proposition}
\newtheorem{Lem}[The]{Lemma}
\theoremstyle{definition}
\newtheorem{Def}[The]{Definition}
\newtheorem{Rem}[The]{Remark}
\numberwithin{equation}{section}
\numberwithin{The}{section}

\newcommand{\be}{\begin{eqnarray}}
\newcommand{\ee}{\end{eqnarray}}

\newcommand{\by}{\begin{eqnarray*}}
\newcommand{\ey}{\end{eqnarray*}}

\newcommand{\bn}{\begin{enumerate}}
\newcommand{\en}{\end{enumerate}}

\newcommand{\bi}{\begin{itemize}}
\newcommand{\ei}{\end{itemize}}

\def\Fbar{{\overline F}}

\def\frac#1#2{{#1 \over #2}}

\def \P {\mathbb{P}}

\def\xx{\boldsymbol{x}}
\def\XX{\boldsymbol{X}}
\def\mmu{\boldsymbol{\mu}}
\def\ddelta{\boldsymbol{\delta}}

\def\ttheta{\boldsymbol{\theta}}
\def\yy{\boldsymbol{y}}
\def\zz{\boldsymbol{z}}
\def\YY{\boldsymbol{Y}}

\def\ww{\boldsymbol{w}}
\def\vv{\boldsymbol{v}}

\def\Zero{\boldsymbol{0}}

\begin{document}
\title{On Rapid Variation of Multivariate Probability Densities}
\author{
Haijun Li\footnote{{\small\texttt{lih@math.wsu.edu}}, Department of Mathematics and Statistics, Washington State University, 
Pullman, WA 99164, U.S.A.}
}
\date{April 2021}
\maketitle

\begin{abstract}
Multivariate rapid variation describes decay rates of joint light tails of a multivariate distribution. We impose a local uniformity condition to control decay variation of distribution tails along different directions, and using higher-order tail dependence of copulas, we prove that a rapidly varying multivariate density implies rapid variation of the joint distribution tails. As a corollary, rapid variation of skew-elliptical distributions is established under the assumption that the underlying density generators belong to the max-domain of attraction of the Gumbel distribution.



\end{abstract}

\section{Introduction and preliminaries}
\label{S1}

Multivariate regular variation has been extensively studied in the literature; see, e.g., \cite{Resnick07}, for detailed discussions and various applications in analyzing multivariate extremes. Regular variation of multivariate heavy tails of a distribution on $\mathbb{R}^d$ enjoys a scaling property, that has proven useful in multivariate asymptotic theory,  and in particular, such a scaling property allows to establish the conditions in terms of multivariate densities which imply
 regular variation of multivariate distribution tails \cite{HO1984, HR1987}. 

The goal of this paper is to establish a similar closure property for multivariate light tails and apply it to studying {\em rapid variation}, in the sense of de Haan \cite{DH1970},  of skew-elliptical distributions. We 
develop a local uniformity condition, under which {rapidly varying} multivariate densities imply the {rapid variation} for multivariate distribution tails. In contrast to the scaling property of multivariate regular variation, multivariate rapid variation satisfies only an additive stability condition (see \cite{JL19}), which lacks strong scaling homogeneity that may be needed for establishing the closure property. To overcome this difficulty, we utilize higher-order tail dependencies of copulas, introduced and studied in \cite{HJ10,HJL12, LH14}, and their relations to multivariate extremes. It is known that (first-order) tail dependence of copulas is equivalent to multivariate regular variation of distributions on $[0,\infty]^d$ with tail equivalent univariate marginal distributions \cite{LS09}, and together with this,  
our copula approach used in this paper captures the universality of multivariate dependence for distributions  with regularly or rapidly varying marginal tails of different types, and further illustrates its usefulness in multivariate analysis.   

Our research is motivated from the study on multivariate tail behaviors of skew-elliptical distributions, reported in \cite{JL19}, where tail densities of skew-elliptical distributions and their associated copulas are derived for both cases of regular and rapid variation. Since most multivariate distributions are specified directly by densities, it is important to have closure criteria in terms of densities which imply regular or rapid variation of multivariate distribution tails. Such closure criteria were obtained in \cite{HO1984, HR1987} for multivariate regular variation, and our results in this paper fill the void for multivariate rapid variation, showing further that skew-elliptical distributions are rapidly varying if their density generators belong to the max-domain of attraction of the Gumbel distribution.

The paper is organized as follows. Section 2 establishes the relation of local uniformity for the joint tails between multivariate rapid variation and its copula, and then proves the closure theorem for a rapidly varying density to ensure rapid variation of the multivariate distribution tails. Section 3 extends a result of \cite{JL19} on tail densities of skew-elliptical distributions with rapidly varying density generators to the entire $\mathbb{R}^d$ and establishes rapid variation of the joint tails of skew-elliptical distributions. Skew-elliptical distributions studied here belong to a general class of selection distributions \cite{BD2001}, and these multivariate distributions can account for both skewness and heavy/light tails, and have found widespread application in various areas \cite{Azzalini2013}. Section 4 concludes the paper with some remarks. In the rest of this section, some notation and preliminaries on higher-order tail dependencies of copulas  are highlighted (see \cite{HJL12, LH14} for details). 

Any vector in $\mathbb{R}^d$ or $\mathbb{R}^d_+ = [0,\infty)^d$ is denoted by row vector $\xx = (x_1, \dots, x_d)$ and its transpose, a column vector, is denoted by $\xx^\top$. Two measurable functions $f(t,x)$ and $g(t,x)$ are tail equivalent, and denoted as $f(t,x)\sim g(t,x)$, if $f(t,x)/g(t,x)\to 1$ as $t\to \infty$ for each $x$. These functions are said to be locally uniformly tail equivalent, and denoted as $f(t,x)\sim_{l.u.} g(t,x)$, if $f(t,x)/g(t,x)\to 1$ as $t\to \infty$ locally uniformly in $x$. 
Operations of vectors (sums, products, etc.) and relations (inequalities, intervals, etc.) are taken component-wise. 
We consider throughout this paper that any involved slowly varying function is continuous. The assumption is rather mild due to 
 Karamata's representation (see, e.g., \cite{Resnick07}) that any slowly  varying function can be written as the product of a continuous function and a measurable function with  positive constant limit.

A copula $C$ is a multivariate distribution with uniformly distributed univariate marginal distributions on $[0,1]$. Sklar's theorem (see, e.g.,  \cite{Joe2014}) states that every multivariate distribution $F$ with univariate marginal distributions $F_1,\dots, F_d$ can be written as $F(x_1,\dots, x_d)=C(F_1(x_1),\dots,F_d(x_d))$ for some $d$-dimensional copula $C$. In fact, in the case of continuous univariate margins, $C$ is unique and 
$$C(u_1,\dots,u_d)={F}(F_1^{-1}(u_1), \dots, F_d^{-1}(u_d)), \ (u_1, \dots, u_d)\in [0,1]^d.$$
Let $(U_1,\dots,U_d)$ denote a random vector with distribution $C$ and $U_i, 1\le i\le d$, being uniformly distributed on $[0,1]$. The survival copula $\widehat{C}$ is defined as follows:
\begin{equation}
	\label{dual}
	\widehat{C}(u_1,\dots,u_d)=\P(1-U_1\le u_1,\dots,1-U_d\le u_d)=\overline{C}(1-u_1,\dots,1-u_d)
\end{equation}
where $\overline{C}$ is the joint survival function of $C$.  Assume throughout this paper that the density $c(\cdot)$ of copula $C$ exists, and that $c(\cdot)$ is continuous in some small open neighborhoods of $\boldsymbol 0$ and ${\boldsymbol 1}=(1, \dots, 1)$. (Its dimensionality is always clear from the context without explicit mention.) 
\begin{Def}
	\label{d2.1}
	\begin{enumerate}
		\item The upper {tail density of $C$ with tail order $\kappa_U$} 
		is defined as follows:
		\begin{equation}
		\lambda_U({\boldsymbol w}; \kappa_U):= \lim_{u\to 0^+}
		\frac{c(1-uw_i, 1\le i\le d)}{u^{\kappa_U-d}\ell(u)}>0, ~ {\boldsymbol w}=(w_1,\dots,w_d)\in [0,\infty)^d\backslash \{\boldsymbol 0\},\label{U tail}
		\end{equation}
		provided that the non-zero limit exists for some $\kappa_U\ge 1$ and for some function $\ell(\cdot)$ that is slowly varying at $0$.
		\item The lower {tail density of $C$ with tail order $\kappa_L$} 
		is defined as follows:
		\begin{equation}
		\lambda_L({\boldsymbol w}; \kappa_L):= \lim_{u\to 0^+}
		\frac{c(uw_i, 1\le i\le d)}{u^{\kappa_L-d}\ell(u)}>0, ~ {\boldsymbol w}=(w_1,\dots,w_d)\in [0,\infty)^d\backslash \{\boldsymbol 0\},\label{L tail}
		\end{equation}
		provided that the non-zero limit exists for some $\kappa_L\ge 1$ and for some function $\ell(\cdot)$ that is slowly varying at $0$.
	\end{enumerate}
\end{Def}

\noindent
Here and hereafter $\ell(u)$ is said to be slowly varying at $0$ if $\ell(u^{-1})$ is slowly varying at $\infty$; that is, $\ell(ux)/\ell(u)\to 1$, $x>0$, as $u\to 0$. 

\begin{Def}
	\label{d2.2}
	\begin{enumerate}
		\item  The upper tail dependence function with tail order $\kappa_U$ is defined as follows:
		\begin{equation*}
			%\label{U tail order}
			b_U(\ww;\kappa_U):=\lim_{u\to 0^+}\frac{\overline{C}(1-uw_i, 1\le i\le d)}{u^{\kappa_U} \ell(u)},~w = (w_1, \dots, w_d)\in \mathbb{R}_+^d, 
		\end{equation*}
		provided that the non-zero limit exists at any continuity point $\ww$ of $b_U(\cdot; \kappa_U)$,  
		for $\kappa_U\ge 1$ and for some function $\ell(u)$ that is slowly varying at $0$.
		\item  The lower tail dependence function with tail order $\kappa_L$ is defined as follows:
		\begin{equation*}
		%\label{L tail order}
		b_L(\ww;\kappa_L):=\lim_{u\to 0^+}\frac{{C}(uw_i, 1\le i\le d)}{u^{\kappa_L} \ell(u)},~w = (w_1, \dots, w_d)\in \mathbb{R}_+^d, 
		\end{equation*}
		provided that the non-zero limit exists at any continuity point $\ww$ of $b_L(\cdot; \kappa_L)$,  
		for $\kappa_L\ge 1$ and for some function $\ell(u)$ that is slowly varying at $0$.
	\end{enumerate}
\end{Def}

The fact that $\kappa_U\ge 1$ or $\kappa_L\ge 1$ follows from the Fr\'echet-Hoeffding upper bound. 
Clearly, $\kappa_U\le d$ (or $\kappa_L\le d$) for copulas with positively  upper (or lower) orthant dependence. In addition, the tail densities and tail order functions are scale-homogeneous with the following scaling properties
\begin{equation}
	\label{scaling1}
	\lambda_U(t{\boldsymbol w}; \kappa_U)=t^{\kappa_U-d}\lambda_U({\boldsymbol w}; \kappa_U), ~~\lambda_L(t{\boldsymbol w}; \kappa_L)=t^{\kappa_L-d}\lambda_L({\boldsymbol w}; \kappa_L), \ {\boldsymbol w}\in [0,\infty)^d\backslash\{\boldsymbol 0\}
\end{equation}
\begin{equation*}
	%\label{scaling2}
	b_U(t{\boldsymbol w}; \kappa_U)=t^{\kappa_U}b_U({\boldsymbol w}; \kappa_U), ~~b_L(t{\boldsymbol w}; \kappa_L)=t^{\kappa_L}b_L({\boldsymbol w}; \kappa_L), \ {\boldsymbol w}\in \mathbb{R}^d_+
\end{equation*}
for any $t>0$. It is seen from these scaling properties that these functions are completely determined by values of these functions at $\ww=(w_1, \dots, w_d)$, $0< w_i\le 1$, $1\le i\le d$.

\begin{Pro}
	\label{p1}\rm 
	Assume that tail dependence functions and tail densities exist. 
	\begin{enumerate}
		\item If  \eqref{U tail} holds locally uniformly at $\ww \in (0,\infty)^d$, then 
		\[\lambda_U(\ww;\tau_U)=\frac{\partial^d b_U(\ww;\tau_U)}{\partial w_1 \cdots \partial w_d},\ \ww=(w_1, \dots, w_d)\in \mathbb{R}_+^d\backslash\{0\}.
		\]
		\item If  \eqref{L tail} holds locally uniformly at $\ww \in (0,\infty)^d$, then 
		\[\lambda_L(\ww;\tau_L)=\frac{\partial^d b_L(\ww;\tau_L)}{\partial w_1 \cdots \partial w_d},\ \ww=(w_1, \dots, w_d)\in \mathbb{R}_+^d\backslash\{0\}.
		\]
	\end{enumerate}
\end{Pro}
The above results are proved in \cite{LH14}. Since the survival copula $\widehat{C}$ (see \eqref{dual}) can be used to transform lower tail properties of $(U_1,\dots,U_d)$ into the corresponding upper tail properties of $(1-U_1,\dots,1-U_d)$, the lower tail property (2) of Proposition \ref{p1} can be obtained immediately from the upper tail case. Because of this duality, we focus on upper tails only in this paper.

\section{Rapid variation for multivariate densities}
\label{S3}

Let $F$ be a $d$-dimensional distribution with density $f$, copula $C$  and continuous marginal distributions $F_1, \dots, F_d$. 
With proper translation and scaling sequences, marginal distributions $F_i$'s belong to the max-domain of attraction of one of the three distribution families; that is, Fr\'echet, Gumbel or Weibull distribution.
We focus on the Gumbel case for light tails in this paper.  Assume  that for the upper tail case, marginal distributions $F_1, \dots, F_d$ are  right-tail equivalent in the sense that
\begin{equation}
\label{equiv}
f_i(t)\sim a_if_1(t),~\mbox{as}~t\to \infty,~~1\le i\le d, 
\end{equation}
where $f_i$ is the density of the $i$-th marginal distribution $F_i$, and $0<a_i<\infty$ is a constant,  $1\le i\le d$,  with $a_1=1$. Note that \eqref{equiv} implies the usual tail equivalence of the marginal survival functions, 
\begin{equation}
	\label{dist equiv}
	\Fbar_i(t)\sim a_i\Fbar_1(t),~\mbox{as}~t\to \infty,~~1\le i\le d, 
\end{equation}
where $\overline{F}_i(t) = 1-F_i(t)$, $1\le i\le d$, but the reverse may not be true in general.

Assume now that the marginal density $f_i$, $1\le i\le d$, is continuous and satisfies
\begin{equation}
	\label{G tail}
	f_i(t+m_i(t)x)\sim f_i(t)e^{-x}\sim a_if_1(t)e^{-x}, ~x\in \mathbb{R},~~\mbox{as}~t\to \infty, 
\end{equation}
where the self-neglecting function $m_i(\cdot)$ can be taken to be a differentiable function for which its  derivative converges to $0$ (\cite{BGT1987, HF2006}). According to Bloom's theorem (see \cite{Bloom76}), 
\begin{equation}
\label{Bloom}
\lim_{t\to \infty}\frac{m_i(t+m_i(t)x)}{m_i(t)}=1, \ 1\le i\le d,
\end{equation}
hold locally uniformly in $x\in \mathbb{R}$. All these functions are known to belong to the gamma class (see \cite{Omey13}). The gamma class $\Gamma_\alpha(m)$ consists of all the measurable functions $g$, denoted by $g\in \Gamma_\alpha(m)$, for which there exists a
measurable and positive function $m$ such that
\[\lim_{t\to \infty}\frac{g(t+m(t)x)}{g(t)}=e^{\alpha x},\ \ x\in \mathbb{R}.
\]
The following two results, obtained in \cite{Omey13}, are used in deriving our main results.
\begin{Lem}\rm 
	\label{Omey1} If $g\in \Gamma_{-\kappa}(m)$, $\kappa>0$, where $m(\cdot)$ is self-neglecting, then  
\[\frac{g(t+m(t)x)}{g(t)}\to e^{-\kappa x}, \kappa>0,\ \mbox{locally uniformly in}\ x\in \mathbb{R}. 
\]
\end{Lem}
\begin{Lem}\rm 
	\label{Omey2} Suppose that $m(\cdot)$ is self-neglecting. Then 
	${g(t+m(t)x)}\sim {g(t)}$
	if and only if $g(x) = \ell(G(x))$, where $\ell(\cdot)$ is slowly varying at $0$ and $G\in \Gamma_{-\kappa}(m)$, $\kappa>0$. 
\end{Lem}
\noindent
It is worth mentioning that Omey in \cite{Omey13} obtained the stronger results, but the descriptions in Lemmas \ref{Omey1} and \ref{Omey2} suit our purpose.  

Lemma \ref{Omey1} implies that the convergences
\begin{equation*}
%\label{Omey}
\lim_{t\to \infty}\frac{f_i(t+m_i(t)x)}{f_i(t)}=e^{-x}, \ 1\le i\le d
\end{equation*}
hold locally uniformly in $x\in \mathbb{R}$. Since the derivative of $m_i(t)$ converges to $0$, the derivative of $t+m_i(t)x$ with respect to $t$ converges to 1, locally uniformly in $x$. Integrating both sides of \eqref{G tail} with respect to $t$ implies that 
\[\Fbar_i(t+m_i(t)x)\sim_{l.u.} \Fbar_i(t)e^{-x}, ~x\in \mathbb{R},~~\mbox{as}~t\to \infty.  
\]
That is, $F_i$ is in the max-domain of attraction of the Gumbel distribution, and thus the reciprocal hazard rate $\Fbar_i(t)/f_i(t)$ can be taken as $m_i(\cdot)$ and $m_i(t)\sim m_1(t)$ as $t\to \infty$, $1\le i\le d$.

\begin{Def}
	\label{d3.1} Suppose that $(X_1, \dots, X_d)$ has a distribution $F$ with density $f$, that is tail equivalent in the sense of \eqref{equiv} and satisfies \eqref{G tail}. 
	\begin{enumerate}
		\item The tail density $\lambda(\cdot)$ at $\infty$ is defined as
		\begin{equation}
		\lambda({\boldsymbol x}):=\lim_{t\to \infty}\frac{f(t+m(t)x_1, \dots, t+m(t)x_d)}{m^{-d}(t)V^{\kappa_U}(t)}, \ \ \xx\in \mathbb{R}^d, 
		\label{tail density def}
		\end{equation}
		where $\kappa_U>0$,  $m(\cdot)$ is self-neglecting and $V\in \Gamma_{-1}(m)$, provided that the non-zero limit exists. 
		\item The distribution $F$ is said to have a rapidly varying tail at $\infty$ if there exists a non-null Radon measure $\mu(\cdot)$ such that for all the $\mu$-continuity points $\xx\in \mathbb{R}^d$	(satisfying that  $\mu\big(\partial (\xx, +\infty]\big)=0$, where $\partial B$ denotes the boundaries of set $B$), 	
		\begin{equation}
		\label{rapidly verying tail}
	\frac{\mathbb{P}\big(X_1>t+m(t)x_1, \dots, X_d>t+m(t)x_d\big)}{V^{\kappa_U}(t)}\to \mu\Big(\prod_{i=1}^d(x_i, +\infty]\Big), \ \mbox{as $t\to \infty$}, 
		\end{equation} 
		where $\kappa_U>0$,  $m(\cdot)$ is self-neglecting and $V\in \Gamma_{-1}(m)$. 
	\end{enumerate}
\end{Def}

\begin{Rem}
	\label{G00}
	\begin{enumerate}
		\item The self-neglecting function $m(t)$ in \eqref{tail density def} and \eqref{rapidly verying tail} can be taken as $m(t) =m_1(t) = \overline{F}_1(t)/f_1(t)$, $t\ge 0$. The function $V(t)$ in \eqref{tail density def} and \eqref{rapidly verying tail} can be taken as $\Fbar_1(t)
		[\ell(\Fbar_1(t))]^{1/\kappa_U}$, for some function  $\ell(\cdot)$ that is slowly varying at $0$. It follows from the Karamata representation for slowly varying functions that 
		\[V(t+m(t)x)=\Fbar_1(t+m(t)x)
		[\ell(\Fbar_1(t+m(t)x))]^{1/\kappa_U}\]
		\[ \sim \Fbar_1(t)e^{-x}
		[\ell(\Fbar_1(t)e^{-x})]^{1/\kappa_U}
		\sim \Fbar_1(t)e^{-x}
		[\ell(\Fbar_1(t))]^{1/\kappa_U}=V(t)e^{-x},
		\]
		and by Lemma \ref{Omey1}, the convergence is locally uniform. Scaling functions, such as $V(t)$, is unique up to a constant. 
		\item  In contrast to the tail density in multivariate regular variation \cite{HR1987}, $\lambda(\cdot)$ in \eqref{tail density def} does not have a scaling property, and it satisfies the following weak additive stability condition:
		\[\lambda(\xx + z{\boldsymbol 1})=\lambda(\xx)e^{-z\kappa_U}, \ {\boldsymbol x}\in \mathbb{R}^d, 
		\]
		for any $z\in \mathbb{R}$. 
	\end{enumerate}
\end{Rem}

The following result is needed in analyzing multivariate tails of $F$ via our copula approach.
\begin{Pro}
	\label{G} 
	\rm 
	 If the limit \eqref{U tail} holds locally uniformly in ${\boldsymbol w}\in (0,\infty)^d$, then $F$ has a tail density $\lambda(\cdot)$  that is related to the upper tail density $\lambda_U(\cdot;\kappa_U)$ of $C$ as follows: for any ${\boldsymbol x}=(x_1, \dots, x_d)\in \mathbb{R}^d$, 
	\begin{eqnarray}
		\lambda({\boldsymbol x})
		&=&\prod_{i=1}^da_i\,e^{-\sum_{i=1}^dx_i}\lambda_U(a_1e^{-x_1}, \dots, a_de^{-x_d}; \kappa_U)\nonumber\\
		&=&\lambda_U(a_1e^{-x_1}, \dots, a_de^{-x_d}; \kappa_U)|J(a_1e^{-x_1}, \dots, a_de^{-x_d})|,\label{tail density transform1}
	\end{eqnarray}
	where  
	$J(a_1e^{-x_1}, \dots, a_de^{-x_d})$ is the Jacobian determinant of the homeomorphic transform $w_i=a_ie^{-x_i}$, $1\le i\le d$.
\end{Pro}

This result was obtained in \cite{LH14} on $[0,\infty)^d\backslash \{\boldsymbol 0\}$, but can be extended to $\mathbb{R}^d$ using the same proof. In fact, a stronger result in which the convergence for a tail density holds locally uniformly on $\mathbb{R}^d$ can be obtained as follows.

\begin{Pro}
	\label{G1} 
	\rm 
	The limit \eqref{U tail} holds locally uniformly in ${\boldsymbol w}\in (0,\infty)^d$ if and only if the limit \eqref{tail density def} holds locally uniformly in ${\boldsymbol x}\in \mathbb{R}^d$. If either of these two conditions holds, then \eqref{tail density transform1} holds. 
\end{Pro}

\noindent
{\sl Proof.} (1) Suppose that  \eqref{U tail} holds locally uniformly in ${\boldsymbol w}\in (0,\infty)^d$, with a slowly varying function $\ell(\cdot)$. Consider, for any ${\boldsymbol x}=(x_1, \dots, x_d)\in \mathbb{R}^d$, that 
\[\frac{f(t+m(t)x_1, \dots, t+m(t)x_d)}{m^{-d}(t)V^{\kappa_U}(t)}
=\frac{c\big(F_1(t+m(t)x_1), \dots, F_d(t+m(t)x_d)\big)}{[\Fbar_1(t)]^{\kappa_U-d}\ell(\Fbar_1(t))}\frac{\prod_{i=1}^df_i(t+m(t)x_i)}{\prod_{i=1}^df_1(t)},
\]
where $m(t)$ and $V(t)$ are described as these in Remark \ref{G00} (1). 
By the tail equivalency \eqref{equiv} and Lemma \ref{Omey1}, we have
\[\frac{f_i(t+m(t)x_i)}{f_1(t)}\to a_ie^{-x_i}, \ \ \frac{\overline{F}_i(t+m(t)x_i)}{\overline{F}_1(t)}\to a_ie^{-x_i}, \ \ 1\le i\le d, 
\]
converge locally uniformly in $\xx\in \mathbb{R}^d$. Since $c(\cdot)$ is ultimately continuous at $\boldsymbol 1$, the Heine-Cantor theorem implies that for any small $\epsilon>0$, there exists an $N_1$ such that when $t>N_1$, for all $\xx\in B\subset \mathbb{R}^d$, where $B$ is compact, 
\[\left\vert\frac{c\Big(1-\frac{\Fbar_1(t+m(t)x_1)}{\Fbar_1(t)}\Fbar_1(t), \dots, 1-\frac{\Fbar_d(t+m(t)x_d)}{\Fbar_1(t)}\Fbar_1(t)\Big)}{[\Fbar_1(t)]^{\kappa_U-d}\ell(\Fbar_1(t))}\frac{\prod_{i=1}^df_i(t+m(t)x_i)}{\prod_{i=1}^df_1(t)}\right. \]
\[\left.
-\frac{c(1-a_ie^{-x_i}\overline{F}_1(t), 1\le i\le d)}{[\Fbar_1(t)]^{\kappa_U-d}\ell(\Fbar_1(t))}\prod_{i=1}^da_ie^{-x_i}\right\vert \le \frac{\epsilon}{2}. 
\]
Since  \eqref{U tail} holds locally uniformly in ${\boldsymbol w}\in (0,\infty)^d$, there exists an $N_2$, such that as $t>N_2$, for all $\xx\in B\subset \mathbb{R}^d$, 
\[\left|\frac{c(1-a_ie^{-x_i}\overline{F}_1(t), 1\le i\le d)}{[\Fbar_1(t)]^{\kappa_U-d}\ell(\Fbar_1(t))}\prod_{i=1}^da_ie^{-x_i}-\lambda_U(a_ie^{-x_i}, 1\le i\le d; \kappa_U)\prod_{i=1}^da_ie^{-x_i}\right|\le \frac{\epsilon}{2}. 
\]
Therefore, as $t>\max\{N_1, N_2\}$, for all $\xx\in B\subset \mathbb{R}^d$, where $B$ is compact, 
\[\left|\frac{f(t+m(t)x_1, \dots, t+m(t)x_d)}{m^{-d}(t)V^{\kappa_U}(t)}-\lambda_U(a_ie^{-x_i}, 1\le i\le d; \kappa_U)\prod_{i=1}^da_ie^{-x_i}\right|\le \epsilon
\]
which implies that the limit \eqref{tail density def} holds locally uniformly in $\xx\in \mathbb{R}^d$ and \eqref{tail density transform1} holds.

(2) Suppose that  \eqref{tail density def} holds locally uniformly in ${\boldsymbol x}\in \mathbb{R}^d$. Let $L(t):= V^{\kappa_U}(t)/\overline{F}_1^{\kappa_U}(t)$, and obviously, $L(t+m(t)x)\sim L(t)$ where $m(t)$ is self-neglecting. By Lemma \ref{Omey2}, $L(t)=\ell(\overline{F}_1(t))$\footnote{The function $G$ in Lemma \ref{Omey2} must be $\overline{F}_1(t)$, by virtue of the construction used in the proof of Lemma \ref{Omey2}.}, where $\ell(\cdot)$ is slowly varying at $0$. Also observe that $a_ie^{-x_i}=w_i$ if and only if $x_i=-\log(w_i/a_i)$ for $-\infty < x_i<\infty$ and $0< w_i< \infty$. Let $u=\overline{F}_1(t)$, and consider, for $0< w_i< \infty$, 
	\begin{eqnarray*}
&&\frac{c(1-uw_i, 1\le i\le d)}{u^{\kappa_U-d}\ell(u)} = \frac{c(1-\overline{F}_1(t)a_ie^{-x_i}, 1\le i\le d)}{m^{\kappa_U-d}(t)f_1^{\kappa_U-d}(t)\ell(\overline{F}_1(t))} \\
&=& \frac{c\Big(1-\frac{\overline{F}_1(t)a_ie^{-x_i}}{\overline{F}_i(t+m(t)x_i)}\overline{F}_i(t+m(t)x_i), 1\le i\le d\Big)}{m^{-d}(t)\Fbar_1^{\kappa_U}(t)\ell(\overline{F}_1(t))}\frac{\prod_{i=1}^df_1 (t)}{\prod_{i=1}^df_i(t+m(t)x_i)}\prod_{i=1}^df_i(t+m(t)x_i). 
\end{eqnarray*}
It follows from the tail equivalency and Lemma \ref{Omey1} that
\[\frac{\overline{F}_1(t)a_ie^{-x_i}}{\overline{F}_i(t+m(t)x_i)}\to 1, \ 1\le i\le d, \ \ \frac{\prod_{i=1}^df_1 (t)}{\prod_{i=1}^df_i(t+m(t)x_i)}\to \prod_{i=1}^d a_i^{-1}e^{x_i}, 
\]
converge locally uniformly in ${\boldsymbol x}\in \mathbb{R}^d$. Since $c(\cdot)$ is ultimately continuous at $\boldsymbol 1$, the Heine-Cantor theorem implies that for any small $\epsilon>0$, there exists an $N_1$ such that when $t>N_1$, for all $\xx\in B\subset \mathbb{R}^d$, where $B$ is compact, 
\[\left|\frac{c(1-uw_i, 1\le i\le d)}{u^{\kappa_U-d}\ell(u))}-
\frac{c\Big(1-\overline{F}_i(t+m(t)x_i), 1\le i\le d\Big)}{m^{-d}(t)\Fbar_1^{\kappa_U}(t)\ell(\overline{F}_1(t))}\prod_{i=1}^df_i(t+m(t)x_i)\prod_{i=1}^d a_i^{-1}e^{x_i}\right|\le \frac{\epsilon}{2}. 
\]
Since $f(t+m(t)x_1, \dots, t+m(t)x_d) = c\big(1-\overline{F}_i(t+m(t)x_i), 1\le i\le d\big)\prod_{i=1}^df_i(t+m(t)x_i)$, it follows from the local uniform convergence of \eqref{tail density def} that there exists an $N_2$ such that as $t>N_2$, for all $\xx\in B\subset \mathbb{R}^d$, where $B$ is compact, 
\[\left|\frac{f(t+m(t)x_1, \dots, t+m(t)x_d)}{m^{-d}(t)V^{\kappa_U}(t)}\prod_{i=1}^d a_i^{-1}e^{x_i}-\lambda(x_1, \dots, x_d)\prod_{i=1}^d a_i^{-1}e^{x_i}\right|\le \frac{\epsilon}{2}, 
\]
where $V(t)=\Fbar_1(t)[\ell(\overline{F}_1(t))]^{1/\kappa_U}$. Therefore, as $t>\max\{N_1, N_2\}$, for all $\xx\in B\subset \mathbb{R}^d$, where $B$ is compact, 
\[\left|\frac{c(1-uw_i, 1\le i\le d)}{u^{\kappa_U-d}\ell(u))}-
\lambda(x_1, \dots, x_d)\prod_{i=1}^d a_i^{-1}e^{x_i}\right|\le {\epsilon}, 
\]
which implies that \eqref{U tail} holds locally uniformly on $({\boldsymbol 0},\infty)$, and \eqref{tail density transform1} holds for $0< w_i = a_ie^{-x_i}< \infty$, $1\le i\le d$. \hfill $\Box$

The conditions of local uniform convergences in Proposition \ref{G1} are rather mild for the functions with multiplicative scaling properties. In fact, the local uniform convergences can be imposed only on a relatively compact subset of $\mathbb{R}^d_+$.  
\begin{Pro}
	\label{G2} 
	\rm 
If the limit \eqref{U tail} holds locally uniformly in ${\boldsymbol w}\in \prod_{i=1}^d(0,a_i]$, $a_i>0$, $1\le i\le d$, then the limit \eqref{U tail} holds locally uniformly in ${\boldsymbol w}\in (0,\infty)^d$. 
\end{Pro} 

\noindent
{\sl Proof:} 
For any $\ww\in (0,\infty)^d$, let $[\ww]:=||\ww||_2/\min\{a_i\}$ for the  Euclidean norm $||\cdot||_2$. Obviously, $0< w_i/[\ww]\le a_i$, $1\le i\le d$. Consider
\[\lim_{u\to 0^+}
\frac{c(1-uw_i, 1\le i\le d)}{u^{\kappa_U-d}\ell(u)}=[\ww]^{\kappa_U-d}\lim_{u\to 0^+}
\frac{c(1-u[\ww]w_i/[\ww], 1\le i\le d)}{(u[\ww])^{\kappa_U-d}\ell(u[\ww])}\frac{\ell (u[\ww])}{\ell (u)}.
\]
For any compact subset $B\subset (0,\infty)^d$, the set $\{\ww/[\ww]: \ww\in B\}$ is a compact subset of $\prod_{i=1}^d(0,a_i]$. 
Since \eqref{U tail} holds locally uniformly on $\prod_{i=1}^d(0,a_i]$, 
\[
\frac{c(1-u[\ww]w_i/[\ww], 1\le i\le d)}{(u[\ww])^{\kappa_U-d}\ell(u[\ww])}
\]
converges locally uniformly in $\ww\in (0,\infty)^d$. In addition, it follows from the local uniform convergence of a slowly varying function that ${\ell (u[\ww])}/{\ell (u)}$ converges to 1 locally uniformly in $\ww\in \mathbb{R}_+^d\backslash \{{\boldsymbol 0}\}$. It then follows from \eqref{scaling1} that 
\[[\ww]^{\kappa_U-d}\lim_{u\to 0^+}
\frac{c(1-u[\ww]w_i/[\ww], 1\le i\le d)}{(u[\ww])^{\kappa_U-d}\ell(u[\ww])}\frac{\ell (u[\ww])}{\ell (u)}=[\ww]^{\kappa_U-d}\lambda_U({\boldsymbol w}/[\ww]; \kappa_U)=\lambda_U({\boldsymbol w}; \kappa_U),
\]
locally uniformly in $\ww\in (0,\infty)^d$.
\hfill $\Box$

The other ingredient of our copula method is the integration form of Proposition \ref{p1}. 

\begin{Pro}
	\label{t1}\rm 
	Assume that tail dependence functions and tail densities of a copula $C$ exist. 
	\begin{enumerate}
		\item If  \eqref{U tail} holds locally uniformly in $\ww\in (0,\infty)^d$, then 
		\begin{equation}
			\label{e11}
			b_U(\ww;\tau_U) = \int_0^{w_1}\cdots \int_0^{w_d}\lambda_U(\xx;\tau_U)\,{\rm d}\xx,\ \ww=(w_1, \dots, w_d)\in \mathbb{R}_+^d.
		\end{equation}
		\item If  \eqref{L tail} holds locally uniformly in $w\in (0,\infty)^d$, then 
		\[b_L(\ww;\tau_L) = \int_0^{w_1}\cdots \int_0^{w_d}\lambda_L(\xx;\tau_L)\,{\rm d}\xx,\ \ww=(w_1, \dots, w_d)\in \mathbb{R}_+^d.
		\]
	\end{enumerate}
\end{Pro}

\noindent
{\sl Proof.} Only (1) is proved and the proof of (2) is similar. Since the limit \eqref{U tail} is locally uniform, Proposition \ref{p1} (1) holds. 

The local uniformity of \eqref{U tail} implies that the tail density $\lambda_U(\cdot;\tau_U)$ is integrable on the compact subset 
$\{\xx\le \ww: \wedge_{i=1}^dx_i\ge \epsilon\}$ for any fixed small $\epsilon>0$ and any fixed $\ww\in \mathbb{R}_+^d\backslash\{0\}$, where $\wedge_{i=1}^dx_i=\min\{x_i, 1\le i\le d\}$. Since \eqref{U tail} converges locally uniformly and $c(\cdot)$ is ultimately continuous at $\boldsymbol 1$, $\partial ^{|S|}b_U(\ww;\tau_U)/\partial \ww_{S}$ is continuous in $\ww$, for any subset $S \subseteq \{1, \dots, d\}$. Because $\mathbb{R}^d$ is locally compact, the Heine-Cantor theorem implies that $\partial ^{|S|}b_U(\ww;\tau_U)/\partial \ww_{S}$ is locally uniformly continuous in $\ww$, for any subset $S \subseteq \{1, \dots, d\}$,  which further implies via differentiablity of local uniform convergence that 
\begin{equation}
	\label{e10}
	\lim_{\vv_{S^c}\to \ww_{S^c}}\frac{\partial^{|S|}b_U(\ww_S,\vv_{S^c}; \tau_U)}{\partial \ww_S}=\frac{\partial^{|S|}}{\partial \ww_S}\lim_{\vv_{S^c}\to \ww_{S^c}}b_U(\ww_S,\vv_{S^c}; \tau_U)=\frac{\partial^{|S|}b_U(\ww_S,\ww_{S^c}; \tau_U)}{\partial \ww_S}
\end{equation}
with the initial condition that $b_U(\ww_S, {\boldsymbol 0}_{S^c}; \tau_U)=0$, for any subset $S \subset \{1, \dots, d\}$.

It follows from Fatou's lemma that
\by
\int_{[{\boldsymbol 0}, \ww]}\lambda_U(\xx;\tau_U)\,{\rm d}\xx&=& \int \lim_{\epsilon\to 0}I_{\{\xx\le \ww, \wedge_{i=1}^dx_i\ge \epsilon\}}\lambda_U (\xx; \tau_U)\,{\rm d}\xx\\
&\le & {\liminf_{\epsilon\to 0}}\int_{\{\xx\le \ww, \wedge_{i=1}^dx_i\ge \epsilon\}}\lambda_U (\xx; \tau_U)\,{\rm d}\xx
\ey
\by
&\le & {\liminf_{\epsilon\to 0}}\lim_{u\to 0}\int_{\{\xx\le \ww, \wedge_{i=1}^dx_i\ge \epsilon\}}\frac{c(1-ux_i, 1\le i\le d)}{u^{\kappa_U-d}\ell(u)}\,{\rm d}\xx\\
&\le & \lim_{u\to 0}\int_{\{{\boldsymbol 0}\le \xx\le \ww\}}\frac{c(1-ux_i, 1\le i\le d)}{u^{\kappa_U-d}\ell(u)}\,{\rm d}\xx\\
&=& \lim_{u\to 0}\frac{\overline{C}(1-uw_i, 1\le i\le d)}{u^{\kappa_U}\ell (u)}=b_U(\ww;\tau_U )<\infty,
\ey
that is, $\lambda_U(\cdot;\tau_U)$ is integrable on $[{\boldsymbol 0}, \ww]$. Taking integrations on the both sides of the expression in Proposition \ref{p1} (1) and \eqref{e10} yield \eqref{e11}. 
\hfill $\Box$

The main result of this section is the following. 

\begin{The}
	\label{pdf CDF} \rm 
	Let $(X_1, \dots, X_d)$ have a distribution $F$ with density $f$ and copula $C$, satisfying marginal tail equivalency \eqref{equiv}.  Suppose that tail density 
	\[\lambda({\boldsymbol w})=\lim_{t\to \infty}\frac{f(t+m(t)w_1, \dots, t+m(t)w_d)}{m^{-d}(t)V^{\kappa_U}(t)}>0, \ \ \ww\in \mathbb{R}^d, 
	\]
	exists, for a self-neglecting function $m(t)$, and a rapidly varying function $V(t)$ and $\kappa_U>0$. 
	If this limit holds locally uniformly in $\ww\in \mathbb{R}^d$, then ${F}$ is rapidly varying at $\infty$, i.e., 
	\[\lim_{t\to \infty}\frac{\mathbb{P}\big(X_1>t+m(t)x_1, \dots, X_d>t+m(t)x_d\big)}{V^{\kappa_U}(t)}=  \int_{[\xx, \infty)}\lambda(\ww)\,{\rm d}\ww, \ \xx\in \mathbb{R}^d.  
	\]
\end{The}

\noindent
{\sl Proof.} It follows from Proposition \ref{G1} that the limit \eqref{U tail} holds locally uniformly in ${\boldsymbol w}\in (0,\infty)^d$ and \eqref{tail density transform1} holds with $V(t) = \overline{F}_1(t) \ell(\overline{F}_1(t))^{1/\kappa_U}$. By Proposition \ref{t1} (1), we have
\[b_U(\ww;\tau_U) = \int_0^{w_1}\cdots \int_0^{w_d}\lambda_U(\yy;\tau_U)\,{\rm d}\yy,\ \ww=(w_1, \dots, w_d)\in \mathbb{R}_+^d.
\]
By setting $y_i=a_ie^{-v_i}$, $1\le i\le d$, observe that, for any $0< w_i< \infty$, $1\le i\le d$,  
\begin{eqnarray*}
b_U(\ww;\tau_U) &=& \int_{-\log (w_1/a_1)}^\infty \cdots \int_{-\log (w_d/a_d)}^\infty \lambda_U (a_1e^{-v_1}, \dots, a_de^{-v_d}; \kappa_U)\prod_{i=1}^da_ie^{-v_i}\,{\rm d}\vv\\
&=& \int_{-\log (w_1/a_1)}^\infty \cdots \int_{-\log (w_d/a_d)}^\infty \lambda (\vv)\,{\rm d}\vv. 
\end{eqnarray*}
Let $x_i=-\log (w_i/a_i)$, $1\le i\le d$, we have
\[b_U(\ww;\tau_U)= \int_{x_1}^\infty \cdots \int_{x_d}^\infty \lambda (\vv)\,{\rm d}\vv, \ \  \xx\in \mathbb{R}^d. 
\]
On the other hand, by letting $u=\overline{F}_1(t)$, we have
\begin{eqnarray*}
b_U(\ww;\tau_U) &=& \lim_{u\to 0^+}\frac{\overline{C}(1-uw_i, 1\le i\le d)}{u^{\kappa_U} \ell(u)}= \lim_{t\to \infty}\frac{\overline{C}(1-\overline{F}_1(t)w_i, 1\le i\le d)}{\overline{F}_1(t)^{\kappa_U} \ell(\overline{F}_1(t))}\\
&=& \lim_{t\to \infty}\frac{\overline{C}(1-\overline{F}_i(t)e^{-x_i}, 1\le i\le d)}{\overline{F}_1(t)^{\kappa_U} \ell(\overline{F}_1(t))}=\lim_{t\to \infty}\frac{\overline{C}(1-\overline{F}_i(t+m(t)x_i), 1\le i\le d)}{\overline{F}_1(t)^{\kappa_U} \ell(\overline{F}_1(t))}\\
&=& \lim_{t\to \infty}\frac{\mathbb{P}\big(X_1>t+m(t)x_1, \dots, X_d>t+m(t)x_d\big)}{V^{\kappa_U}(t)}, 
\end{eqnarray*}
where the third and fourth equalities follow from the tail equivalence, the monotonicity and continuity of $\overline{C}(\cdot)$.   The integration formula then follows.  
\hfill $\Box$

\begin{Rem}
	\begin{enumerate}
		\item Since $1-F(\xx)$ is non-increasing in $\xx$, it follows from Dini's theorem and Moore-Osgood theorem that the multivariate regular varying distribution $F$ studied in \cite{HR1987} clearly satisfies the tail equivalency \eqref{dist equiv} for its marginal survival functions. In contrast, Theorem \ref{pdf CDF} on multivariate rapid variation assumes a stronger tail equivalency \eqref{equiv} for univariate marginal densities. 
		\item The multivariate regular variation studied in \cite{HR1987} describes the upper tail dependence within the cone $(0,\infty)^d$. The multivariate rapid variation \eqref{tail density def} can be used to analyzing higher order upper tail dependence (see \cite{HJL12, LH14}), and our closure theorem on multivariate rapid variation describes tail decay rates on $[\xx, \infty)$, that contains, with an appropriate $\xx\in \mathbb{R}^d$, any relatively compact subset of the entire $[-\infty, \infty]^d\backslash \{-\infty\}$. 
	\end{enumerate}
\end{Rem}

\section{Rapid variation for skew-elliptical distributions}
\label{S4}

A general class of skew-elliptical distributions can be constructed   using conditioning from symmetrical elliptical distributions (see \cite{BD2001} for details).   A $d$-dimensional random vector ${  \XX}$ is said to be {elliptically distributed}  with mean vector ${  \mmu}\in \mathbb{R}^d$ and nonnegative-definite dispersion matrix $\Sigma\in \mathbb{R}^{d\times d}$, if the characteristic function $\varphi_{{  {\XX-\mmu}}}$ of $  {\XX-\mmu}$ is a function of the quadratic form ${ \zz}\Sigma{\,  \zz}^\top$, where, in this paper, $\zz^\top$ denotes the column vector, that is the transpose of a row vector $\zz\in \mathbb{R}^d$.
That is, $\varphi_{  {\XX-\mmu}}({ \zz})=\phi({ \zz}\Sigma{\, \zz}^\top)$,  $\zz\in \mathbb{R}^d$,  for a measurable function $\phi: [0,\infty)\to \mathbb{R}$, which is called a  characteristic generator. We assume that the density of $\XX$ exists and is given, for all $\xx\in \mathbb{R}^d$, by
$$
f_{\XX}(\xx; \mmu, \Sigma)= |\Sigma|^{-1/2}g_{d}\big((\xx-\mmu)\Sigma^{-1}(\xx-\mmu)^\top\big),
$$
where $g_{d}: [0,\infty)\to \mathbb{R}_+$, is known as the $d$-dimensional density generator that satisfies 
\begin{equation}
	\int_0^\infty r^{d/2-1}g_{d}(r)\,{\rm d}r = {\Gamma(d/2)}/{\pi^{d/2}}. 
	\label{generator integral}
\end{equation}
The elliptically distributed random vectors with densities  are denoted as ${  \XX}\sim_d E_d(\mmu, \Sigma, g_{d})$. A comprehensive review of  characterizations and  properties of elliptical distributions can be found in \cite{FKN90}.  We assume throughout that the density generator $g_{d}$ is continuous over $(0, \infty)$ and that $g_d (t)$ is eventually non-increasing as $t\to \infty$. Note that this assumption is rather weak, due to \eqref{generator integral},  and all the commonly used density generators are continuous functions that eventually {decrease} to zero.

Consider now a $(d+1)$-dimensional random vector $(X_0, \XX)\sim_d E_{d+1}({  \mmu}^*, \Sigma^*, g_{d+1})$, where $\XX=(X_1, \dots, X_d)$,
$\mmu^*=(0, \mmu)$, $\mmu=(\mu_1, \dots, \mu_d)$,  $g_{d+1}$ is a $(d+1)$-dimensional density generator, and the dispersion matrix $\Sigma^*$ has the form
\[
\Sigma^*=
\begin{pmatrix}
	1 & \ddelta \\
	\ddelta^\top & \Sigma
\end{pmatrix}
\]
with $\ddelta = (\delta_1, \dots, \delta_d)\in \mathbb{R}^d$,  such that $\Sigma$ is positive-definite and $\ddelta \Sigma^{-1}\ddelta^\top<1$ \cite{BD2001}. The distribution of $\YY = [\XX|X_0>0]$ is called in \cite{BD2001} a skew-elliptical distribution and denoted as $\YY\sim_d SE_{d}({  \mmu}, \Sigma, g_{d+1}, \ddelta)$, where $\ddelta$ is known as the vector of skewness parameters. Note that marginally, $\XX\sim_d E_d({  \mmu}, \Sigma, g_{d})$, where, for all $s \in [0, \infty)$, 
\begin{equation}
	\label{g}
	g_{d}(s)= 
	2\int_0^\infty g_{d+1}(r^2+s)\,{\rm d}r;  
\end{equation}
see Eq.~(2.23) in \cite{FKN90}.
It follows from Theorem 2.18 of \cite{FKN90} that  $\mathbb{P}(X_0>0|\XX=\yy)=F\big((\yy-\mmu)\ttheta^\top; g_{q(\yy)}\big)$, where $\ttheta= {\ddelta\Sigma^{-1}}/{(1-\ddelta \Sigma^{-1}\ddelta^\top)^{1/2}}$ and $F(\cdot\ ; g_{q(\yy)})$ is the cumulative distribution function of the univariate elliptical distribution $E_1(0, 1, g_{q(\yy)})$ with $q(\yy)=(\yy-\mmu)\Sigma^{-1}(\yy-\mmu)^\top$ and $g_{q(\yy)}(s)= {g_{d+1}\big(s+q(\yy)\big)}\big/{g_{d}\big(q(\yy)\big)}$, 
for all $s \in [0, \infty)$.  If $\ddelta = \Zero$, then  $\ttheta=\Zero$ and thus $f_{\YY}=f_{\XX}(\cdot; \mmu, \Sigma)$ in this case. In general, $f_{\YY}$ is skewed. The density of $\YY$ can be written, for all $\yy\in \mathbb{R}^d$, as
\begin{eqnarray}
	\label{density Y}
	f_{\YY}(\yy)&=&2f_{\XX}(\yy; \mmu, \Sigma)F\big((\yy-\mmu)\ttheta^\top; g_{q(\yy)}\big)\nonumber\\
		\label{density Y explicit}
	&=&2|\Sigma|^{-1/2}\int_{-\infty}^{(\yy-\mmu)\ttheta^\top}g_{d+1}\big(r^2+(\yy-\mmu)\Sigma^{-1}(\yy-\mmu)^\top \big)\,{\rm d}r. 
\end{eqnarray}
{Examples} of \eqref{density Y} include skew-normal and skew-$t$ distributions \cite{BD2001}.

\begin{Def}
	\label{GMD Gumbel}
	A measurable function $g:[0,\infty)\to \mathbb{R}_+$ is said to be in the $d$-dimensional, quadratic max-domain of attraction for the Gumbel distribution with auxiliary function $m$  if for any $d\times d$ nonnegative-definite matrix $Q$ and for any $\xx\in \mathbb{R}^d$, 
	\begin{equation}
		\label{Gumbel generator}
		g \big( \{ t{\boldsymbol 1}+m(t)\xx\} Q \{ t{\boldsymbol 1}+m(t)\xx\}^\top \big)\sim 	g \big(t^2{\boldsymbol 1}Q{\boldsymbol 1}^\top \big)
		\exp(-\xx Q {\boldsymbol 1}^\top)
	\end{equation}
	as $t \to \infty$, for some self-neglecting function $m$, i.e., $m \big( t+m(t)x\big) \sim m(t)$, $x\in \mathbb{R}$, as $t\to \infty$. 
\end{Def}

\begin{Rem}
	\label{GMD Gumbel remark}
	This notion is slightly more general than that introduced in \cite{JL19}, by requiring only nonnegative definiteness on $\mathbb{R}^d$. Because of this, 
	if $g$ is in the $d$-dimensional, quadratic max-domain of attraction for the Gumbel distribution with auxiliary function $m$, then \eqref{Gumbel generator} clearly also holds for any $l\times l$ nonnegative-definite matrix $Q$ and $\xx \in \mathbb{R}^l$, where $1\le l\le d-1$. 
\end{Rem}

\begin{Pro}
	\label{P3-1} 
	\rm If a continuous function $g$ is in the $d$-dimensional, quadratic max-domain of attraction for the Gumbel distribution with auxiliary function $m$, then
	 	\begin{equation}
	 	\label{Gumbel generator l.u.}
	 	\frac{g \big( \{ t{\boldsymbol 1}+m(t)\xx\} Q \{ t{\boldsymbol 1}+m(t)\xx\}^\top \big)}{g (t^2{\boldsymbol 1}Q{\boldsymbol 1}^\top )}\to 	
	 	\exp(-\xx Q {\boldsymbol 1}^\top), 
	 \end{equation}
 locally uniformly in $\xx\in \mathbb{R}^d$, for any $d\times d$ nonnegative-definite matrix $Q$. 
\end{Pro}

\noindent
{\sl Proof.} Let 
\[H(t,\xx) = \frac{g\big(t^2\{{\boldsymbol 1}+m(t)\xx/t\}Q\{{\boldsymbol 1}+m(t)\xx/t\}^\top\big)}{g(t^2{\boldsymbol 1}Q{\boldsymbol 1}^\top)}, \ t\ge 0, \xx\in \mathbb{R}^d. 
\]
Because of \eqref{Gumbel generator}, for each $\zz\in B$, where $B\subset \mathbb{R}^d$ is  compact, $H(t,\zz)$ is Cauchy, $t\to \infty$. For any small $\epsilon>0$, there exists a large $N_{\zz}$, such that whenever $t, t'>N_{\zz}$, 
\begin{equation}
	\label{cauchy1}
\big|H(t,\zz)-H(t',\zz)\big|\le \epsilon/4. 
\end{equation}
Since $g$ and $m$ are continuous, the function $H(t,\xx)$  
is continuous on $[0,\infty)\times \mathbb{R}^d$. For the same $\epsilon>0$, $\zz\in B$ and $t'>N_{\zz}$, since $m(t)/t\to 0$ as $t\to \infty$, choose a small $\delta_{\zz} > 0$, such that whenever $||\zz-\xx||_2\le \delta_{\zz}$ for the Euclidean norm $||\cdot||_2$, we also have
\begin{equation}
	\label{cauchy2}
	\big|H(t',\zz)-H(t',\xx)\big|\le \epsilon/4. 
\end{equation}
Combining \eqref{cauchy1} and \eqref{cauchy2}, we have that 
whenever $t, t'> N_{\zz}$ and $||\zz-\xx||_2\le \delta_{\zz}$,  
\begin{equation}
	\label{cauchy3}
\big|H(t,\zz)-H(t',\xx)\big|\le \epsilon/2.
\end{equation}
Note that the collection of balls $B(\zz, \delta_{\zz}/2)$ covers $B$, so by compactness there is a finite subcover, say involving $\zz_1, \dots, \zz_n$, such that $\cup_{i=1}^n B(\zz_i, \delta_{\zz_i}/2)\supset B$. 

Let $\delta = \min_{1\le i\le n}\delta_{\zz_i}/2$ and $N=\max_{1\le i\le n}N_{\zz_i}$. Suppose that $||\xx-\yy||_2\le \delta$ for arbitrary $\xx, \yy\in B$ and $t, t'> N$. By the choice of $\zz_1, \dots, \zz_n$ and the triangle inequality, there exists an $i$ 
such that $||\xx-\zz_i||_2\le \delta_{\zz_i}$ and $||\zz_i-\yy||_2\le \delta_{\zz_i}$. It follows from \eqref{cauchy3} that 
\[\big|H(t,\xx)-H(t',\yy)\big| \le \big|H(t,\xx)-H(t,\zz_i)\big|+\big|H(t,\zz_i)-H(t',\yy)\big|\le \frac{\epsilon}{2}+\frac{\epsilon}{2}=\epsilon. 
\]
Hence, $H(t,\xx)$ is uniformly Cauchy for all $\xx\in B$, as $t\to \infty$, where $B\subset \mathbb{R}^d$ is  compact. Therefore, \eqref{Gumbel generator l.u.} converges locally uniformly in $\xx\in \mathbb{R}^d$. 
\hfill $\Box$

\begin{Rem}
	\begin{enumerate}
		\item Bloom \cite{Bloom76} proved in the univariate case  \eqref{Bloom} holds locally uniformly if $m$ is continuous, and one can use a similar idea to prove Proposition \ref{P3-1} for a multivariate case. But our proof of Proposition \ref{P3-1} is direct and elementary. 
		\item Omey \cite{Omey13} obtained the local uniformity in Lemma \ref{Omey1} for the univariate case. It is not clear whether local uniformity holds in the multivariate case \eqref{Gumbel generator l.u.} without the continuity assumption. 
	\end{enumerate}
\end{Rem}

\begin{Pro}
	\label{P3-2} 
	\rm If a $(d+1)$-dimensional density generator $g_{d+1}$ is in the quadratic max-domain of attraction for the Gumbel distribution with auxiliary function $m$,  then the $d$-dimensional density generator $g_{d}$, defined by \eqref{g}, is also in the quadratic max-domain of attraction for the Gumbel distribution with same auxiliary function $m$. 
\end{Pro}

\noindent
{\sl Proof.} For any positive-definite matrix $\Sigma\in \mathbb{R}^{d\times d}$, it follows from Remark \eqref{GMD Gumbel remark} and Proposition \ref{P3-1} that
\begin{equation}
	\label{g converge 0}
	\frac{g_{d+1} \big(\{ t{\boldsymbol 1}+m(t)\ww \} \Sigma^{-1} \{ t{\boldsymbol 1}+m(t)\ww \} ^{\top}\big)}{g_{d+1} \big(t{\boldsymbol 1} \Sigma^{-1} t{\boldsymbol 1}^{\top}\big)}\to \exp(-\ww \Sigma^{-1} {\boldsymbol 1}^\top), 
\end{equation}
locally uniformly in $\ww\in \mathbb{R}^d$.
Define
\[Q_s= 
\begin{pmatrix}
	s^2 & {\boldsymbol 0}\\
	{\boldsymbol 0}^\top & \Sigma^{-1}
\end{pmatrix}. 
\]
For any $s\in \mathbb{R}$, $Q_s$ is nonnegative-definite. Therefore, by Proposition \ref{P3-1}, 
\[\frac{g_{d+1} \big( \{ t+m(t)\times 0, t{\boldsymbol 1}+m(t)\ww\}
Q_s \{ t+m(t)\times 0, t{\boldsymbol 1}+m(t)\ww\}^{\top}\big)}{g_{d+1} \big(t^2(1,{\boldsymbol 1})Q_s(1, {\boldsymbol 1})^\top \big)} \to \exp\{- (0,\ww) Q_s(1, {\boldsymbol 1})^\top\}
\]
locally uniformly in $\ww\in \mathbb{R}^d$. Observe that $\exp\{-(0,\ww) Q_s(1, {\boldsymbol 1})^\top\} = \exp (-\ww\Sigma^{-1}{\boldsymbol 1}^\top)$, and thus
\[\frac{g_{d+1} \big(
	t^2s^2+ \{ t{\boldsymbol 1}+m(t)\ww \} \Sigma^{-1} \{ t{\boldsymbol 1}+m(t)\ww \} ^{\top}\big)}{g_{d+1} \big(
	t^2s^2+ t{\boldsymbol 1} \Sigma^{-1} t{\boldsymbol 1}^{\top}\big)}
\]
\begin{equation}
	\label{g converge}
=\frac{g_{d+1} \big(
	t^2s^2+ t\{ {\boldsymbol 1}+m(t)\ww/t \} \Sigma^{-1} t\{ {\boldsymbol 1}+m(t)\ww/t \} ^{\top}\big)}{g_{d+1} \big(
	t^2s^2+ t{\boldsymbol 1} \Sigma^{-1} t{\boldsymbol 1}^{\top}\big)}\to \exp(-\ww \Sigma^{-1} {\boldsymbol 1}^\top), 
\end{equation}
locally uniformly in $\ww\in \mathbb{R}^d$.

Observe that the quadratic forms $t\{ {\boldsymbol 1}+m(t)\ww/t \} \Sigma^{-1} t\{ {\boldsymbol 1}+m(t)\ww/t \} ^{\top}$ and $t{\boldsymbol 1} \Sigma^{-1} t{\boldsymbol 1}^{\top}$ are strictly positive and go to positive infinity, as $t\to \infty$. Also observe that $g_{d+1}$ and $m$ are locally uniformly continuous by the Heine-Cantor theorem. 
Since $m(t)/t\to 0$  and $t^2s^2\ge 0$ for any $s\in \mathbb{R}$, \eqref{g converge 0} implies that \eqref{g converge} holds for all $s\in \mathbb{R}$ uniformly. 
 That is, for $\ww\in B\subset \mathbb{R}$, where $B$ is compact, and for any small $\epsilon>0$, there exists a large positive $N(B, \epsilon)$ such that whenever $t>N(B, \epsilon)$, 
\[(1-\epsilon) g_{d+1} \big(
t^2s^2+ t{\boldsymbol 1} \Sigma^{-1} t{\boldsymbol 1}^{\top}\big)\exp(-\ww \Sigma^{-1} {\boldsymbol 1}^\top)\le 
g_{d+1} \big(
t^2s^2+ \{ t{\boldsymbol 1}+m(t)\ww \} \Sigma^{-1} \{ t{\boldsymbol 1}+m(t)\ww \} ^{\top}\big)  
\]
\begin{equation}
	\label{g converge key}
\le (1+\epsilon) g_{d+1} \big(
t^2s^2+ t{\boldsymbol 1} \Sigma^{-1} t{\boldsymbol 1}^{\top}\big)\exp(-\ww \Sigma^{-1} {\boldsymbol 1}^\top), \ \forall\ s\in \mathbb{R}. 
\end{equation}
Taking the integrals with respect to $s$ from $0$ to $\infty$ and then letting $t\to \infty$ and $\epsilon\to 0$, we obtain that 
\[\frac{\int_0^\infty g_{d+1} \big(
	t^2s^2+ \{ t{\boldsymbol 1}+m(t)\ww \} \Sigma^{-1} \{ t{\boldsymbol 1}+m(t)\ww \} ^{\top}]\big)\, {\rm d}s}{\int_0^\infty g_{d+1} \big(
	t^2s^2+ t{\boldsymbol 1} \Sigma^{-1} t{\boldsymbol 1}^{\top}\big) \, {\rm d}s}\to \exp(-\ww \Sigma^{-1} {\boldsymbol 1}^\top), 
\]
locally uniformly in $\ww\in \mathbb{R}^d$. Applying  \eqref{g} leads to 
\[\frac{g_d\big(\{ t{\boldsymbol 1}+m(t)\ww \} \Sigma^{-1} \{ t{\boldsymbol 1}+m(t)\ww \} ^{\top}\big)}{g_d\big(t{\boldsymbol 1} \Sigma^{-1} t{\boldsymbol 1}^{\top}\big)}\to \exp(-\ww \Sigma^{-1} {\boldsymbol 1}^\top), 
\]
locally uniformly in $\ww\in \mathbb{R}^d$, and the result follows.
\hfill $\Box$

\begin{The}
	\label{skewed CDF-1} \rm 
	Let $(Y_1, \dots, Y_d)$ have a skew-elliptical distribution $F$ with density $f$ given by \eqref{density Y explicit}.  If the density generator $g_{d+1}$ is in the quadratic max-domain of attraction for the Gumbel distribution with auxiliary function $m$,  then the tail density $\lambda(\ww)$ at $\infty$ is given by
	\begin{equation}
		\lambda({\boldsymbol w})=\lim_{t\to \infty}\frac{f(t{\boldsymbol 1}+m(t)\ww)}{m^{-d}(t)V^{\kappa_U}(t)}=	\left\{\begin{array}{ll}
			2|\Sigma|^{-1/2}\exp(-\ww\Sigma^{-1}{\boldsymbol 1}^\top) & \mbox{if ${\boldsymbol 1}\ttheta^\top\ne {0}$,}\\
			|\Sigma|^{-1/2}\exp(-\ww\Sigma^{-1}{\boldsymbol 1}^\top) & \mbox{if ${\boldsymbol 1}\ttheta^\top={0}$,}
		\end{array}
		\right. \label{eq-lambda Gumbel}
	\end{equation}
where the convergence holds locally uniformly 
in  $\ww\in \mathbb{R}^d$, and $\kappa_U={\boldsymbol 1}\Sigma^{-1}{\boldsymbol 1}^\top>0$. 
\end{The}

\noindent
{\sl Proof.} Assume, without loss of generality, that $\mmu=\Zero$. 
Consider, from \eqref{density Y explicit}, that for any fixed $\ww\in \mathbb{R}^d$,
\begin{align}
	f_{\YY} (t{\boldsymbol 1}+m(t)\ww )&= 2|\Sigma|^{-1/2}\int_{-\infty}^{\{t{\boldsymbol 1}+m(t)\ww\}\ttheta^{\top}}
	g_{d+1} \big( r^2+ \{ t{\boldsymbol 1}+m(t)\ww \} \Sigma^{-1} \{ t{\boldsymbol 1}+m(t)\ww\} ^{\top}\big)\,
	{\rm d}r\nonumber\\
	&= 2|\Sigma|^{-1/2}G(t,\ww), 
	\label{f_Y Gumbel}
\end{align}
where 
\begin{eqnarray*}
G(t,\ww)&=& \int_{-\infty}^{\{t{\boldsymbol 1}+m(t)\ww\}\ttheta^{\top}}g_{d+1} \big(
r^2+ \{ t{\boldsymbol 1}+m(t)\ww \} \Sigma^{-1} \{ t{\boldsymbol 1}+m(t)\ww \} ^{\top}\big)\,{\rm d}r\\
&=& t\int_{-\infty}^{\{{\boldsymbol 1}+m(t)\ww/t\}\ttheta^{\top}}g_{d+1} \big(
t^2s^2+ \{ t{\boldsymbol 1}+m(t)\ww \} \Sigma^{-1} \{ t{\boldsymbol 1}+m(t)\ww \} ^{\top}\big)\,{\rm d}s.
\end{eqnarray*}
Since \eqref{g converge key} holds, we 
take the integrals with respect to $s$ from $-\infty$ to $\{{\boldsymbol 1}+m(t)\ww/t\}\ttheta^{\top}$ and then let $t\to \infty$ and $\epsilon\to 0$, and obtain from \eqref{g converge key} that 
\begin{equation}
	\label{Gumbel generator integral l.u.}
	\frac{\int_{-\infty}^{\{{\boldsymbol 1}+m(t)\ww/t\}\ttheta^{\top}}g_{d+1} \big(
		t^2s^2+ \{ t{\boldsymbol 1}+m(t)\ww \} \Sigma^{-1} \{ t{\boldsymbol 1}+m(t)\ww \} ^{\top}\big)\,{\rm d}s}{\int_{-\infty}^{\{{\boldsymbol 1}+m(t)\ww/t\}\ttheta^{\top}}g_{d+1} \big(
		t^2s^2+ t{\boldsymbol 1} \Sigma^{-1} t{\boldsymbol 1}^{\top}\big)\,{\rm d}s}\to \exp(-\ww \Sigma^{-1} {\boldsymbol 1}^\top), 
\end{equation}
locally uniformly in $\ww\in \mathbb{R}^d$. In other words,
\begin{eqnarray}
	\label{lu+1}
G(t,\ww) &\sim_{l.u.}& \exp(-\ww \Sigma^{-1} {\boldsymbol 1}^\top)\,t\int_{-\infty}^{\{{\boldsymbol 1}+m(t)\ww/t\}\ttheta^{\top}}g_{d+1} \big(
t^2s^2+ t{\boldsymbol 1} \Sigma^{-1} t{\boldsymbol 1}^{\top}\big)\,{\rm d}s,\\
&=& \exp(-\ww \Sigma^{-1} {\boldsymbol 1}^\top)\int_{-\infty}^{\{t{\boldsymbol 1}+m(t)\ww\}\ttheta^{\top}}g_{d+1} \big(
r^2+ t{\boldsymbol 1} \Sigma^{-1} t{\boldsymbol 1}^{\top}\big)\,{\rm d}r. \label{lu-0}
\end{eqnarray}

Consider the two cases. 
\begin{enumerate}
	\item If ${\boldsymbol 1}\ttheta^{\top}\ne 0$, then 
	$$\{t{\boldsymbol 1}+m(t)\ww\}\ttheta^{\top} = t\Big\{{\boldsymbol 1}+\frac{m(t)}{t}\ww \Big\}\ttheta^{\top}\sim_{l.u.} t{\boldsymbol 1}\ttheta^{\top}, $$
	which implies that
	\begin{equation}
		\label{lu-1}
	\int_{-\infty}^{\{{\boldsymbol 1}+m(t)\ww/t\}\ttheta^{\top}}g_{d+1} \big(
	t^2s^2+ t{\boldsymbol 1} \Sigma^{-1} t{\boldsymbol 1}^{\top}\big)\,{\rm d}s\sim_{l.u.}  \int_{-\infty}^{{\boldsymbol 1}\ttheta^{\top}}g_{d+1} \big(
	t^2s^2+ t{\boldsymbol 1} \Sigma^{-1} t{\boldsymbol 1}^{\top}\big)\,{\rm d}s. 
	\end{equation}
	Let $G(t) = \int_{-\infty}^{t{\boldsymbol 1}\ttheta^{\top}}g_{d+1} \big(
	r^2+ t{\boldsymbol 1} \Sigma^{-1} t{\boldsymbol 1}^{\top}\big)\,{\rm d}r$. It follows from \eqref{Gumbel generator integral l.u.} that for any $x\in \mathbb{R}$, 
	\begin{eqnarray*}	
	G(t+m(t)x) &=&  \int_{-\infty}^{\{t+m(t)x\}{\boldsymbol 1}\ttheta^{\top}}g_{d+1} \big(
	r^2+ \{t+m(t)x\}{\boldsymbol 1} \Sigma^{-1} \{t+m(t)x\}{\boldsymbol 1}^{\top}\big)\,{\rm d}r\\
	&=& t\int_{-\infty}^{\{1+m(t)x/t\}{\boldsymbol 1}\ttheta^{\top}}g_{d+1} \big(
	t^2s^2+ \{t+m(t)x\}{\boldsymbol 1} \Sigma^{-1} \{t+m(t)x\}{\boldsymbol 1}^{\top}\big)\,{\rm d}s\\
	&\sim & \exp(-x {\boldsymbol 1}\Sigma^{-1} {\boldsymbol 1}^\top)\,t\int_{-\infty}^{\{1+m(t)x/t\}{\boldsymbol 1}\ttheta^{\top}}g_{d+1} \big(
	t^2s^2+ t{\boldsymbol 1} \Sigma^{-1} t{\boldsymbol 1}^{\top}\big)\,{\rm d}s\\
	&\sim & \exp(-x {\boldsymbol 1}\Sigma^{-1} {\boldsymbol 1}^\top)\,t\int_{-\infty}^{{\boldsymbol 1}\ttheta^{\top}}g_{d+1} \big(
	t^2s^2+ t{\boldsymbol 1} \Sigma^{-1} t{\boldsymbol 1}^{\top}\big)\,{\rm d}s\\
	&\sim & \exp(-x {\boldsymbol 1}\Sigma^{-1} {\boldsymbol 1}^\top)\int_{-\infty}^{t{\boldsymbol 1}\ttheta^{\top}}g_{d+1} \big(
	r^2+ t{\boldsymbol 1} \Sigma^{-1} t{\boldsymbol 1}^{\top}\big)\,{\rm d}r\\
	&=& G(t)\exp(-x {\boldsymbol 1}\Sigma^{-1} {\boldsymbol 1}^\top). 
	\end{eqnarray*}
	Let $V(t) = m^{d/\kappa_U}(t)G^{1/\kappa_U}(t)$, $t\ge 0$, where $\kappa_U={\boldsymbol 1}\Sigma^{-1}{\boldsymbol 1}^\top>0$, and we then have that $V(t+m(t)x) \sim V(t)e^{-x}$, $x\in \mathbb{R}$. Therefore, combining \eqref{f_Y Gumbel}, \eqref{lu-0} and \eqref{lu-1}, we have
	\[\lambda(\ww) = \lim_{t\to \infty}\frac{f_{\YY} (t{\boldsymbol 1}+m(t)\ww )}{m^{-d}(t)V^{\kappa_U}(t)}= 2|\Sigma|^{-1/2}\exp(-\ww \Sigma^{-1} {\boldsymbol 1}^\top),
	\]
	and the convergence holds locally uniformly in $\ww\in \mathbb{R}^d$. 
	\item Suppose that ${\boldsymbol 1}\ttheta^\top=0$. In this case, it follows from \eqref{lu+1} and \eqref{g} that 
	\begin{eqnarray}	
		G(t,\ww) &\sim_{l.u.}& \exp(-\ww \Sigma^{-1} {\boldsymbol 1}^\top)\,t\int_{-\infty}^{\frac{m(t)}{t}\ww \ttheta^\top}g_{d+1} \big(
		t^2s^2+ t{\boldsymbol 1} \Sigma^{-1} t{\boldsymbol 1}^{\top}\big)\,{\rm d}s\nonumber \\
		&\sim_{l.u.}& \exp(-\ww \Sigma^{-1} {\boldsymbol 1}^\top)\,t\int_{-\infty}^{0}g_{d+1} \big(
		t^2s^2+ t{\boldsymbol 1} \Sigma^{-1} t{\boldsymbol 1}^{\top}\big)\,{\rm d}s\nonumber \\
		&\sim_{l.u.}& \exp(-\ww \Sigma^{-1} {\boldsymbol 1}^\top)\int_{-\infty}^{0}g_{d+1} \big(
		r^2+ t{\boldsymbol 1} \Sigma^{-1} t{\boldsymbol 1}^{\top}\big)\,{\rm d}r\nonumber\\
			&=&\frac{1}{2}\exp(-\ww \Sigma^{-1} {\boldsymbol 1}^\top) g_{d}(t{\boldsymbol 1} \Sigma^{-1} t{\boldsymbol 1}^{\top}). 	\label{lu-3}
\end{eqnarray}
	Let $G(t) = g_{d}(t{\boldsymbol 1} \Sigma^{-1} t{\boldsymbol 1}^{\top})$. Because of \eqref{Gumbel generator} and Proposition \ref{P3-2}, this function satisfies that $G(t+m(t)x) \sim G(t)\exp(-x {\boldsymbol 1}\Sigma^{-1} {\boldsymbol 1}^\top)$ for all $x\in \mathbb{R}$. Let $V(t) = m^{d/\kappa_U}(t)G^{1/\kappa_U}(t)$, $t\ge 0$, where $\kappa_U={\boldsymbol 1}\Sigma^{-1}{\boldsymbol 1}^\top>0$, and we have that $V(t+m(t)x) \sim V(t)e^{-x}$, $x\in \mathbb{R}$. Therefore, combining \eqref{f_Y Gumbel}, \eqref{lu-0} and \eqref{lu-3}, we have
	\[\lambda(\ww) = \lim_{t\to \infty}\frac{f_{\YY} (t{\boldsymbol 1}+m(t)\ww )}{m^{-d}(t)V^{\kappa_U}(t)}= |\Sigma|^{-1/2}\exp(-\ww \Sigma^{-1} {\boldsymbol 1}^\top),
	\]
	and the convergence holds locally uniformly in $\ww\in \mathbb{R}^d$. 
	\hfill $\Box$
\end{enumerate}

\begin{Rem}
	\label{R3-10}
	\begin{enumerate}
		\item The scaling function $V(t) = m^{d/\kappa_U}(t)G^{1/\kappa_U}(t)$  can be explicitly obtained from our construction. Since scaling is unique up to a constant, the two cases in \eqref{eq-lambda Gumbel} can be combined into one expression. 
		\item Theorem \ref{skewed CDF-1} extends the result obtained in \cite{JL19} for skew-elliptical densities with rapid varying tails to the entire $\mathbb{R}^d$, where the convergence \eqref{eq-lambda Gumbel} holds locally uniformly under a slightly weaker condition on density generators.
	\end{enumerate}
\end{Rem}

 To show such a skew-elliptical distribution $F$ is rapidly varying at $\infty$, we need to impose the conditions on the vector $\ddelta$ of skewness parameters, to ensure right-tail equivalence. Let $\bar\theta_i=\delta_i/(1-\delta_i^2)^{1/2}$ for all $i \in \{ 1, \ldots, d\}$. 
\begin{itemize}
	\item [(i)]
	If $\bar\theta_i\ge 0$ for all $i \in \{ 1, \ldots, d\}$, then $f_1, \dots, f_d$ are right-tail equivalent in the sense of \eqref{equiv}. 
	\item [(ii)]
	If $\bar\theta_1= \dots = \bar\theta_d< 0$, then $f_1, \dots, f_d$ are right-tail equivalent in the sense of \eqref{equiv}. 
	\end{itemize}
The same conditions (i) and (ii) are used in \cite{JL19} in deriving higher order tail densities of skew-elliptical copulas.  
Theorems \ref{pdf CDF} and \ref{skewed CDF-1} then lead to the main result of this section. 

\begin{The}
	\label{skewed CDF} \rm 
	Let $(Y_1, \dots, Y_d)$ have a skew-elliptical distribution $F$ with density $f$ given by \eqref{density Y explicit}, satisfying right-tail equivalence described in (i) or (ii) above.  If the density generator $g_{d+1}$ is in the quadratic max-domain of attraction for the Gumbel distribution with auxiliary function $m$,  then ${F}$ is rapidly varying at $\infty$, i.e., 
	\[\lim_{t\to \infty}\frac{\mathbb{P}\big(Y_1>t+m(t)x_1, \dots, Y_d>t+m(t)x_d\big)}{V^{\kappa_U}(t)}=  \int_{[\xx, \infty)}\lambda(\ww)\,{\rm d}\ww, \ \xx\in \mathbb{R}^d,  
	\] 
	where $\kappa_U={\boldsymbol 1}\Sigma^{-1}{\boldsymbol 1}^\top>0$ and for $\ww\in \mathbb{R}^d$, 
	\begin{equation*}
		\lambda({\boldsymbol w})=	\left\{\begin{array}{ll}
			2|\Sigma|^{-1/2}\exp(-\ww\Sigma^{-1}{\boldsymbol 1}^\top) & \mbox{if ${\boldsymbol 1}\ttheta^\top\ne {0}$,}\\
			|\Sigma|^{-1/2}\exp(-\ww\Sigma^{-1}{\boldsymbol 1}^\top) & \mbox{if ${\boldsymbol 1}\ttheta^\top={0}$. }
		\end{array}
		\right.
	\end{equation*}
\end{The}

\begin{Rem}
	\label{r3-6}
 If there are some $\bar{\theta}_i$'s with different signs, then marginal densities may not be right-tail equivalent. For example, it follows from \eqref{density Y explicit} with $\mmu = {\boldsymbol 0}$ that 
		\begin{equation*}
		%	\label{Gumbel tail 4}
			f_{Y_i}(t) = 2\int_{-\infty}^{t\bar\theta_i}g_{2}(r^2+t^2)\,{\rm d}r, \ i\in \{1, \dots, d\}. 
		\end{equation*}
	Then right-tail equivalence is violated for the skew-normal distributions having $\bar{\theta}_i$'s with different signs. 
\end{Rem}

\begin{Exa}\rm 
	\label{e3-1}
	Let $\XX=(Z_0, Z_1, Z_2)$ be a random vector with $3$-dimensional normal distribution having a vector of means $\mmu^*=(0, 0, 0)$ with positive-definite correlation matrix
	\[
	\Sigma^* =
	\begin{pmatrix}
		1 & \ddelta  \\
		\ddelta^\top & \Sigma
	\end{pmatrix},
	\]
	where $\ddelta = (\delta_1,\delta_2)$ and for any $\rho \in (-1,1)$, 
	\begin{equation*}
	%	\label{bi-normal}
		\Sigma = 
		\begin{pmatrix}
			1 & \rho \\
			\rho & 1
		\end{pmatrix}, \quad 
		\Sigma^{-1}=\frac{1}{1-\rho^2}
		\begin{pmatrix}
			1 & -\rho \\
			-\rho & 1
		\end{pmatrix}.
	\end{equation*}
The density generator is given by $g_{3}(x)=(2\pi)^{-3/2}\exp(-x/2)$, which  is in the quadratic max-domain of attraction for the Gumbel distribution with auxiliary function $m(t) = t^{-1}$, $t>0$. 

The $i$th univariate marginal density of $Y_i$ is given by
\[f_{Y_i}(y) = 2\phi(y)\Phi(\bar\theta_iy),\quad \bar\theta_i=\delta_i/(1-\delta_i^2)^{1/2},\ i = 1, 2, 
\]
where $\phi$ denotes the standard normal density and $\Phi$ denotes the standard normal cumulative distribution function. 
The joint density of $(Y_1, Y_2)$ is given by
\begin{equation*}
	%\label{skew-normal exa}
	f_2(y_1, y_2)=2\phi_2(y_1, y_2; \Sigma)\Phi \Big(\sum_{i=1}^2\theta_iy_i \Big)
\end{equation*}
with $(\theta_1, \theta_2) = {\ddelta\Sigma^{-1}}{/} {(1-\ddelta \Sigma^{-1}\ddelta^\top)^{1/2}}$, 
where $\phi_2(\cdot; \Sigma)$ is the $2$-dimensional normal density with zero mean and correlation matrix $\Sigma$. If $(\bar\theta_1, \bar\theta_2)\ge (0,0)$, then the marginal densities are right-tail equivalent in the sense of \eqref{equiv}, with $a_1=1$ and, 
\begin{equation*}
%	\label{a_2}
	a_2=\lim_{y\to \infty} \{ f_{Y_2}(y)/f_{Y_1}(y)\} = 
	\left\{\begin{array}{cl}
		1 & \mbox{if $\bar\theta_1=\bar\theta_2=0$ or $\min (\bar\theta_1, \bar\theta_2) >0$,}\\
		1/2 & \mbox{if $\bar\theta_1>0$ and $\bar\theta_2=0$,}\\
		2 & \mbox{if $\bar\theta_1=0$ and $\bar\theta_2>0$.}
	\end{array}
	\right.
\end{equation*}

Consider the case that $(\delta_1, \delta_2)\ge (0,0)$ and $\rho\ge 0$, then the tail density at $\infty$ can be written explicitly as 
\begin{equation*}
%	\label{bi normal tail density}
	\lambda(w_1, w_2)=
	\left\{\begin{array}{ll}
		2(1-\rho^2)^{-1/2}\exp\{-(w_1+w_2)/(1+\rho)\} & \mbox{if $\theta_1>0$ or $\theta_2>0$,}\\
		(1-\rho^2)^{-1/2}\exp\{-(w_1+w_2)/(1+\rho)\} & \mbox{if $\theta_1=\theta_2=0$,}
	\end{array}
	\right.
\end{equation*}
where the upper tail order $\kappa_U=(1,1)\Sigma^{-1}(1,1)^\top=2/(1+\rho)$ (see \cite{JL19} for detail). According to Theorem \ref{skewed CDF}, the distribution of $\YY$ is rapidly varying at $\infty$, i.e., 
\[\lim_{t\to \infty}\frac{\mathbb{P}\big(Y_1>t+m(t)x_1, Y_2>t+m(t)x_2\big)}{V^{\kappa_U}(t)}=  \int_{[\xx, \infty)}\lambda(w_1,w_2)\,{\rm d}w_1{\rm d}w_2, \ \xx\in \mathbb{R}^2,  
\] 
\begin{equation*}
%	\label{bi normal tail CDF}
=	\left\{\begin{array}{ll}
		2(1+\rho)^2(1-\rho^2)^{-1/2}\exp\{-(x_1+x_2)/(1+\rho)\} & \mbox{if $\theta_1>0$ or $\theta_2>0$,}\\
		(1+\rho)^2(1-\rho^2)^{-1/2}\exp\{-(x_1+x_2)/(1+\rho)\} & \mbox{if $\theta_1=\theta_2=0$,}
	\end{array}
	\right.
\end{equation*}
where the scaling function $V(t)$ can be also explicitly obtained using the constructions in the proof of Theorem \ref{skewed CDF}. 
\end{Exa}

\section{Concluding Remarks}

Multivariate regularly varying densities imply multivariate regular variation of the underlying distributions under a local uniformity condition, and as mentioned in \cite{HR1987}, such a condition is used to control variations across different rays. The closure property of preserving regular variation of multivariate distributions from that of joint densities is  useful in analyzing heavy tail phenomena for multivariate extremes \cite{BE07}. In this paper, we obtain the closure property  from joint densities to multivariate distributions that preserves rapid variation in the sense of de Haan \cite{DH1970}, under the local uniformity that controls light tail variations along different directions. We utilize the copula method to extract higher-order scaling properties from weak additive stability that is possessed by multivariate light tails. Analyzing rapidly varying, light distribution tails for multivariate extremes remains topics of our research in the near future. 

Tail dependence parameters have been obtained in \cite{HL02, Schmidt02} for elliptical distributions and in \cite{FS10, Padoan2011} for skew-$t$ distributions. 
Using the de Haan-Resnick closure theorem in \cite{HR1987}, one can obtain general regularly varying distribution limits for skew-elliptical densities with regular varying generators. In contrast to the heavy tail case, we obtain general light tail distribution limits for skew-elliptical densities with rapidly  varying generators, and our result covers  symmetric, elliptical distributions as a special case, which appears to be new in the literature on elliptical distributions.

%\bibliographystyle{apalike}
%\bibliography{tdens-refs} 

\end{document}